\documentstyle[12pt]{amsart}
\input amssym.def
\input amssym.tex
\begin{document}
\def\B{{\cal B}}
\def\O{{\cal O}}
\def\P{{\Bbb P}}
\def\Q{{\Bbb Q}}
\def\C{{\Bbb C}}
\def\Z{{\Bbb Z}}
\def\qed{\hfill\vbox{\hrule\hbox{\vrule\kern3pt\vbox{\kern6pt}\kern3pt\vrule}\hrule}\bigskip}
\newcommand\Bt{{\B_{\mathrm t}}}
\newcommand\Be{{\B_{\mathrm e}}}
\newcommand\Proj{\mathrm{Proj}}
\newcommand\Spec{\mathrm{Spec}}


\theoremstyle{definition}
\newtheorem{prop}{\bf Proposition}[section]
\newtheorem{tma}[prop]{Theorem}
\newtheorem{lma}[prop]{Lemma}
\newtheorem{cor}[prop]{Corollary}
\newtheorem{df}[prop]{Definition}
\newtheorem{MAIN}[prop]{THEOREM}
\newtheorem{rem}[prop]{Remark}
\newtheorem{ex}[prop]{Example}

\title{Low codimension Fano--Enriques threefolds}
\author{Jorge Caravantes}
\email{\tt jorge\_caravante@@mat.ucm.es}
\maketitle

\centerline{\bf Introduction}

In the 1970s, Reid introduced the graded rings method for the explicit classification of surfaces, which he used to produce a list of 95 K3 quasi-smooth hypersurfaces in weighted projective spaces (which were proved to be the only ones). Later, Fletcher used this method to create more lists of different weighted complete intersections. From the K3 surfaces he developed two lists of anticanonically polarised Fano threefolds that have the K3s as hyperplane section. These two lists for Fano threefolds of codimension one and two can be found in {\bf [F]}. Later on, Alt{\i}nok has developed in {\bf [A]} a formula to calculate the Hilbert series of a Fano threefold (which is very important for the graded rings method) and has written a list of codimension three K3 surfaces (which produces a list of codimension three Fano threefolds). All lists are also in {\bf [GRDW]}.

In this paper we deal with Fano--Enriques threefolds (Fano threefolds with a torsion divisor $\sigma$). These are quotients of Fano threefolds under an action by a $\Z/(r)$ group, where $r$ is the order of $\sigma$. We use the above lists of codimension 1, 2 and 3 to give in this paper all possible Fano--Enriques quotients that can be obtained from these lists.

To find these quotients, we  extend Reid's graded rings method and we  find restrictions for the covers. After that, we  test all members of the lists and, for all those members that satisfy the restrictions, we  calculate a Hilbet series to apply the extended graded rings method and search for a quotient.

The distribution of this paper is as follows. In a first section, we give some preliminaries of Fano--Enriques threefolds. We  describe Alt{\i}nok's method to compute the Hibert series for the anticanonical ring of a Fano threefold and the graded rings method in section 2. In section 3 we give an extension of these methods (Alt{\i}nok's and graded rings) to Fano--Enriques threefolds. Finally, in section 4, we  see the complete way to obtain the lists of Fano--Enriques quotiens and give these lists.

Most of this work in this paper was developed during a stay, supported by Maria Sklodowska-Curie graduate school Threefolds in algebraic geometry (3-fAG) reference number MCFH-1999 00687 grants, at the University of Warwick, which provided all computational resources needed. The author wants to thank Gavin Brown, who introduced me to the K3 database and lots of computational possibilities for geometry, and Miles Reid, who proposed me the Fano--Enriques problem, for all their unvaluable help and advise before, during and after such stay. Currently, the author is enjoying a scholarship from Spanish M.E.C.D. supervised by Enrique Arrondo, who contributed to the final presentation of this paper.

\vskip 10pt
\section{Preliminaries}

Throughout this paper, we  work over the complex field. We recall that a singularity of type ${1\over r}(a_1,...,a_l)$ is a point with an analytic neighbourhood which is isomorphic to a neighbourhood of the origin in $\C^l$ under an action by $\Z/(r)$ consisting on multiplying by $(\epsilon^{a_1},...,\epsilon^{a_l})$ where $\epsilon=e^{2\pi i\over r}$. Sometimes we  just say singularity of type ${1\over r}$ without specifying $a_1,...,a_l$. Our varieties will be irreducible algebraic sets with at most isolated cyclic singularities of index ${1\over r}(1,a,-a)$ with $r$ and $a$ coprime. The necessary background about the weighted projective space and quasi-smooth subvarieties can be found in {\bf[D]} or {\bf [F]}.

\begin{df}
A threefold $X$ is {\it Fano} if it has at most isolated quitient singularities of index ${1\over r}(1,a,-a)$, with $r$ and $a$ coprime, and the anticanonical class $-K_X$ is ample.
\end{df}

\begin{df}
A Fano threefold $X$ is {\it Fano--Enriques} if it has a torsion divisor $\sigma$, that is, there exists $r \in \Z^+$ so that $r\sigma \buildrel lin\over \sim 0$
\end{df}

We begin now the standard construction that relates every Fano--Enriques threefold with a Fano threefold.
Let $X$ be a Fano--Enriques threefold $X$ with torsion divisor $\sigma$ of order $r$. It is well known that 
we can define a covering $\pi:Y \to X$ where  
$$Y=\Spec(\O_X\oplus\O_X(\sigma)\oplus \dots \oplus\O_X((r-1)\sigma)).$$

This covering is $r:1$ at each point in which $\sigma$ is a Cartier divisor. Moreover, there is a regular action by $\Z/(r)$ on $Y$ such that $X$ is the quotient of $Y$ by this action. This action can be described this way: For $j \in \Z/(r)$, consider $\epsilon=e^{2\pi i\over r}$. Then:
$$\O_X\oplus\O_X(\sigma)\oplus \dots \oplus\O_X((r-1)\sigma) \buildrel j \over \longrightarrow \O_X\oplus\O_X(\sigma)\oplus \dots \oplus\O_X((r-1)\sigma))$$
where
$$(a_0,\dots,a_{r-1}) \mapsto (a_0, \epsilon^j a_1,\dots,\epsilon^{(r-1)j}a_{r-1}).$$
Then, for $I \in Y=\Spec(\O_X\oplus\O_X(\sigma)\oplus \dots \oplus\O_X((r-1)\sigma))$ its image is the prime ideal $j^{-1}(I)$.

\begin{rem}\label{obs-accion}
If $X$ is quotient of $Y$ by a $\Z/(r)$ action, that action can be extended to the weighted projective space that contains $Y$. The reason is that the representation of $\Z/(r)$ has a reflection in the anticanonical ring of $Y$ (via the action on the sheaf of structure). The action on $R(Y,-K_Y)$ preserves the degrees of the elements. This means that we have an action for every $R(Y,-K_Y)_d=H^0(Y,-dK_X)$. Therefore, each space will be generated by eigenvectors (this is because the automorphism is an $r$-th root of the identity), which means that the matrix of the automorphism of the ring is diagonal after a suitable change of coordinates. Now we  use the Orbifold Riemann-Roch to classify by eigenvalues ($r$-th roots of 1) the generators of $R(Y,-K_Y)$ (which are eigenvectors after the change of coordinates). This  describes the action that, together with $Y$,  determines $X$. 
\end{rem}

\vskip 10pt
\begin{prop}
Let $X$ be a Fano--Enriques threefold with torsion divisor $\sigma$ of order $r$. Let $Y=\Spec(\O_X\oplus\O_X(\sigma)\oplus \dots \oplus\O_X((r-1)\sigma))$ as above. Then: 
\begin{enumerate}
\item All the singularities of $Y$ are terminal.

\item The cover $Y$ is a $\Q$-Fano threefold.

\end {enumerate}
\end{prop}
\begin{pf}

Let $U\subset X$ be a sufficiently small analytic neighbourhood around a singularity $Q$ of type ${1\over r_Q}(b_Q,1,-1)$. We know that Pic($U$) is isomorphic to $\Z/(r_Q)$ and is generated by the restriction to $U$ of the canonical class, so the local expression of the torsion divisor $\sigma_{|U}$ is $l_Q{K_X}_{|U}$ for some $l_Q$ in $\Z/(r_Q)$. Then, $\alpha:={r_Q\over {\mathrm gcd}(r_Q,l_Q)}$ is the order of $\sigma_{|U}$ in Pic($U$) (i.e. $\alpha\sigma$ is Cartier in $Q$) and divides $r$, the order of $\sigma$ in Pic($X$). Let $Y_{\alpha}$ be the cover of $X$ associated to $\alpha\sigma$. Clearly, $Y$ is a cover for $Y_\alpha$ (there is an obvious monomorphism of rings which defines an epimorphism of schemes with a commutative diagram). Therefore, $Q$ comes from ${r\over \alpha}$ different points with disjoint analytic neighbourhoods $V_i$, $i=1,...,{r\over\alpha}$. Observe now that the canonical cover of a ${1\over r_Q}(b_Q,1,-1)$ singularity is an analytic neighbourhood of $O$ in $\C^3$ defined as 
$$\Spec(\O_U\oplus\O_U({K_X}_{|U})\oplus \dots \oplus\O_U((r-1){K_X}_{|U})).$$
Hence, working locally, we easily deduce that each $V_i$ is of type ${1\over {r_Q\over\alpha}}(b_Q,1,-1)$ which proves (1).

Since $\pi$ is an analytically local isomorphism away from the singularities of $X$ (which are isolated), the anticanonical bundle of $Y$ is $\pi^*(-K_X)$ (i.e. $\pi^{-1}(-K_X)\oplus\pi^{-1}(-K_X+\sigma)\oplus\dots\oplus\pi^{-1}(-K_X+(r-1)\sigma)$). Therefore $Y$ is a Fano threefold and we have (2).

\end{pf}

\begin{ex}
 Consider the complete intersection $Y_{2,2,2} \subset \P^6$ with $\Z/(2)$ acting as the multiplication by $+,+,+,+,-,-,-$ (i.e. the generator of $\Z/(2)$ takes $(x_0,x_1,x_2,x_3,x_4,x_5,x_6)$ to $(x_0,x_1,x_2,x_3,-x_4,-x_5,-x_6)$). Then, since a general $Y$ has eight fixed points (the points in the intersection with $\{ x_4=x_5=x_6=0\}$), $X$ has eight singularities of index ${1\over2}(1,1,1)$.
\end{ex}

The goal of this paper is to reverse this process, i.e. to describe the possible $X$'s for a fixed $Y$ with a prescribed $\Z/(r)$ action.

\vskip 10pt
\section{The graded rings method for Fano threefolds}

In this section we  show how to describe Fano threefolds by finding an appropriate ambient space and equations. We  use the graded rings method introduced by Reid in the 1970s. The aim is to manage some numerical data (i.e. the type of the singularities and the selfintersection of the canonical divisor $-K_Y^3\in \Q$) to search for a  Fano threefold embedded in a weighted projective space  by the ample anticanonical divisor. For this purpose we use the Hilbert series associated to the anticanonical divisor (which depends only on the above numerical data) to guess a possible combination of generators and relations for the ring associated to our threefold embedded in the appropriate weighted projective space. In general, for an ample divisor $D$ (for a Fano threefold, we  use $D=-K_X$) we can construct the ring
$$R(Y,D):= \bigoplus_{n \ge 0} H^0(Y,nD).$$

If we can find generators and relations for this ring, then we have a good description of $Y=\Proj(R(Y,D))$.

We  use Hilbert series to find this $R(Y,D)$. The Hilbert series of $Y$ is:
$$P_Y(t)=\sum_{n\ge 0}h^0(Y,nD)t^n.$$

To calculate this series, we can check Orbifold-RR formula from {\bf [R]}:
\begin{multline}
\chi(\O_Y(D))=\chi(\O_Y)+{1\over12}D(D-K_Y)(2D-K_Y)+\\ 
+{1\over12}Dc_2(Y)+\sum_{Q \in \B}c_Q(D)
\end{multline}
where $\B$ is the basket of terminal quotient singularities ${1\over r_Q}(1,a_Q,-a_Q)={1\over r_Q}(b_Q,1,-1)$ (where $a_Qb_Q\equiv 1$ mod $r_Q$) and the contributions are:
$$c_Q(D)=-i_Q{r_Q^2-1 \over 12r_Q}+ \sum_{j=1}^{i_Q-1}{[b_Qj]_{r_Q}(r_Q-[b_Qj]_{r_Q}) \over 2r_Q}$$
where:
\begin{itemize}
\item $i_Q$ is defined by the condition $\O(D)\simeq \O(i_QK_Y)$ near the singularity $Q$.
\item $[a]_r$ is the minimal nonnegative residue of $a$ mod $r$.
\end{itemize}

\begin{rem} 
Since $-K_X$ is ample, $\chi(-nK_Y)=h^0(Y,-nK_Y)$ for $n \ge 0$ by Kodaira vanishing theorem.
\end{rem}

\begin{tma}\label{altinok-tma}
(Alt\i nok) The Hilbert series $\sum_{n\ge 0}h^0(Y,-nK_Y)t^n$ of a Fano threefold $Y$ can be computed as a rational function on $t$ by the formula:
\begin{multline}\label{altinok}
P_Y(t)={1+t\over (1-t)^2}+{t+t^2\over (1-t)^4}{-K_Y^3\over 2}- \\ - \sum_{Q \in \B}{1\over (1-t)(1-t^r)} \sum_{i=1}^{r-1}{[b_Qi]_{r_Q}(r-[b_i]_{r_Q})\over2r}t^i.
\end{multline}
We call {\it numerical data} of a Fano threefold to $\B$ and the selfintersection of the canonical class $-K_Y^3$.
\end{tma}

\begin{rem}\label{selfint}
We recall from {\bf [A]} that the selfintersection of the canonical class of a Fano threefold X is
$$-K_Y^3 = \sum {b_Q(r_Q-b_Q)\over r_Q} +2k$$
for some integer $k$.
\end{rem}

\begin{ex}\label{fano} Now we  illustrate the graded rings method. Consider a hypotheticalano threefold with just a  ${1\over 2}(1,1,1)$ singularity and $-K_Y^3 = {5\over 2}$. Using the Al{\i}nok's formula (\ref{altinok}), we obtain the Hilbert series:
$$P_Y(t):=1+4t+11t^2+24t^3+46t^4+79t^5+126t^6+189t^7+271t^8+374t^9+\dots$$
By definition, the coefficient of $t^d$ is the dimension of the $\C$- vector space of all homogeneous elements in $R(Y,-K_Y)$ of degree $d$. In particular, we have generators $x_1,x_2,x_3,x_4$ in degree one. Since a generator contributes to this series multiplying by ${1\over 1-t^{a}}$, where $a$ is the degree of the generator, we can multiply by $(1-t)^4$ to simplify the series and discover new generators and relations:
$$(1-t)^4P_Y(t)=1+t^2+t^4-t^5+t^6-t^7+t^8-t^9+\dots $$
This shows that there is a new generator $y$ in degree two. Then, we multiply by $(1-t^2)$ and get $1-t^5$.
Therefore, one expects to have a relation in degree five, which means having a hypersurface $Y_5\subset \P(1,1,1,1,2)$. And in fact, any quasi-smooth equation  works and give the numerical data we started from (i.e. a Fano threefold with a singularity of type ${1\over 2}(1,1,1)$ and $-K_Y^3={5\over 2}$).
\end{ex}

\begin{rem}\label{nonspecial}
The expression ``one expects" means, as in {\bf[ABR]}, ``if there are no extra generators and relations". If there is an extra relation among the monomials (products of the generators), we  need a new generator of the same degree to fill the dimension given by Orbifold-RR and the Hilbert series does not change. We  avoid these special cases (see also Remark \ref{ej2}).
\end{rem}

\vskip 10pt
\section{The graded rings method for Fano--Enriques threefolds}

In this section, we  extend the graded rings method to Fano--Enriques threefolds. We  start with an analogue of Alt{\i}nok's formula. We need an useful remark which  explains why formula (\ref{torsion}) has been motivated.

From Remark \ref{obs-accion}, we deduce that the way to generalize Alt{\i}nok's formula is defining a new Hilbert series considering two degrees: the standard one in $\Z$ and other in $\Z/(r)=\{r-{\mathrm{th \ roots \ of \ }} 1\}$:

Clearly, $H^0(Y, -nK_Y)=\bigoplus_{i=0}^{r-1} H^0(X,-nK_X+i\sigma)$. Therefore, $P_Y(t)= \sum_{i=0}^{r-1} P_X^i(t)$ where $P_X^i(t)=\sum_{n\ge 0} h^0(X,-nK_X+i\sigma)t^n$. 

So we can define the new Hilbert series in $\Z[[t,e]]/(e^r-1)$ as $\sum_{i=0}^{r-1} e^iP_X^i(t)$. The main point of this is that the product in the power series ring agrees with the $\Z/(r)$ action (i.e. a generator in $H^0(X,-nK_X+i\sigma)$ contributes to the Hilbert series multiplying by ${1\over (1-e^it^n)} =1+e^it^n+e^{[2i]_r}t^{2n}+...$). We are grading $R(Y,-K_Y)$ in a new way:

\begin{df}
We  define the {\it bidegree} of an element in $H^0(X,-nK_X+i\sigma)\subset R(Y,-K_Y)$ as $(n,i)\in \Z\oplus\Z/(r)$.
\end{df}

Now we need a formula for $\sum h^0(X,-nK_X+\tau)$, where $\tau$ is a numerically trivial divisor (in particular, it can be a torsion divisor). 

\begin{lma}\label{torsion-lma}
Let $X$ be a Fano threefold with a numerically trivial divisor $\tau$ and a basket of singularities $\B$. For every singularity $Q\in \B$, define $l_Q$ such that, locally in $Q$, $\tau\simeq\O(l_QK_X)$. Then we have that $\sum_{n\ge0}h^0(X,-nK_X+\tau)t^n=$
\begin{multline}\label{torsion}
P_X(t)+ \sum_{Q\in\B} \Big( {r_Q-l_Q\over1-t} {r_Q^2-1\over12r_Q} + \\+{1\over1-t^{r_Q}}\sum_{j=0}^{r_Q-1} \big(\sum_{i=j+1}^{j+r_Q-l_Q}{[b_Qi](r_Q-[b_Qi])\over2r_Q}\big)t^j\Big).
\end{multline}
\end{lma} 

\begin{pf}
By Orbifold Riemann-Roch:
$$\chi(-nK_X)=\chi(\O_X)+{2n^3+3n^2+n\over12}(-K_X^3)+{n\over12}(-K_X)c_2+ \sum_{Q\in\B}c_Q(n)$$
where
$$c_Q(n)=-[-n]{r_Q^2-1\over12r_Q} + \sum_{j=1}^{[-n]-1}{[b_Qj](r_Q-[b_Qj])\over 2r_Q}$$
and 
$$\chi(-nK_X+\tau)= \chi(\O_X)+{2n^3+3n^2+n\over12}(-K_X)+ {n\over12}(-K_X)c_2+ \sum_{Q\in\B} c_Q(-nK_X+\tau)$$
where
$$c_Q(-nK_X+\tau)=-[l_Q-n]{r_Q^2-1\over12r_Q} + \sum_{j=1}^{[l_Q-n]-1}{[b_Qj](r_Q-[b_Qj])\over2r_Q}.$$

It is clear that $\chi(-nK_X+\tau)-\chi(-nK_X)$ is the sum in $\B$ of $c_Q(-nK_X+\tau)-c_Q(n)$. This takes the value:
$$\sum_{Q\in\B}\Big( (r_Q-l_Q){r_Q^2-1\over12r_Q}-\sum_{i=n+1}^{n+r_Q-l_Q}{[b_Qj](r_Q-[b_Qj])\over2r_Q}\Big).$$
This expression is clearly periodic with period $r_Q$, so we get (\ref{torsion})
\end{pf}

\begin{rem}
Lemma \ref{torsion-lma} shows that the numerical data of a Fano--Enriques threefold consist of $-K_X^3$, the order of the torsion divisor $r$, and the basket $\B$ divided in $\Bt$ (singularities where the torsion divisor is not trivial) and $\Be$ (rest of the basket). Moreover, for every singularity $Q\in \Bt$ we add the number $l_Q$, which is determined by the local value of the torsion divisor in the singularity $Q$.
\end{rem}

The generalisation of the graded rings method to Fano--Enriques threefolds  comes by applying Lemmas \ref{altinok-tma} and \ref{torsion-lma}. We illustrate the method with an example.

\begin{ex}\label{ej1}
Consider a Fano--Enriques variety $X$ with basket $\B=\{{1\over10}(1,3,7),2\times{1\over 5}(1,2,3)\}$ with respective $l_Q=6,1,1$ and $-K_X^3={1\over2}$. We look for generators as in Example \ref{fano} but paying attention also to the new weight in $\Z/(5)$. From formulas (\ref{altinok}) and (\ref{torsion}), we get the Hilbert series:
$$(1+t+2t^2+5t^3+9t^4+16t^5+25t^6+38t^7+54t^8+74t^9+...) +$$ 
$$+e(t+2t^2+5t^3+9t^4+15t^5+26t^6+38t^7+54t^8+75t^9+...) +$$ 
$$+e^2(3t^2+5t^3+9t^4+16t^5+25t^6+38t^7+54t^8+75t^9+...) +$$
$$+e^3(t+2t^2+5t^3+9t^4+16t^5+25t^6+37t^7+55t^8+75t^9+...) +$$ 
$$+e^4(t+2t^2+5t^3+9t^4+15t^5+26t^6+38t^7+54t^8+75t^9+...)$$

The coefficient of $e^it^j$ is the dimension of the subspace of elements of bidegree $(i,j)$ in $R(X,-K_X,\sigma)$, so we have generators $x_1,x_2,x_3,x_4$ in the subspaces $H^0(\O_X(-K_X))$, $H^0(\O_X(-K_X+\sigma))$, 
$H^0(\O_X(-K_X+3\sigma))$, $H^0(\O_X(-K_X+4\sigma))$ respectively. They have respective bidegrees (1,0), (1,1), (1,3), (1,4). 

A generator in $H^0(\O_X(-jK_X+i\sigma))$ contributes now to the Hilbert series multiplying by 
$$1+e^it^j+e^{[2i]_r}t^{2j}+\dots$$
This is an inverse for $1-e^it^j$ in $\Z[[e,t]]/(e^r-1)$. Thus we multiply in our case by $1-t, 1-et,1-e^3t,1-e^4t$, so we get
$$(1+t^4+t^8+...)+$$
$$+e(t^6+...)+$$
$$+e^2(t^2+-t^7+...)+$$
$$+e^3(t^8+...)+$$
$$+e^4(t^4-t^9+..)$$
This shows that we need a generator $y$ in $H^0(\O_X(-2K_X+2\sigma))$ (of bidegree (2,2)). We now multiply by $1-e^2t^2$ and get:
$$1-t^5,$$
so we have a relation in $H^0(\O_X(-5K_X))$ (i.e. a weighted-homogeneous polynomial of degree 5 which is invariant by the action) . The expression $x_1^5+x_2^5+x_3^5+x_4^5-x_2y^2$ gives a quasi-smooth relation for $Y_5\subset\P(1,1,1,1,2)$, with the action by $\Z/(5)$ consisting on multiplying each coordinate by $(1,\epsilon,\epsilon^3,\epsilon^4,\epsilon^2)$. 

Now we should check the singularities. The only fixed points by the action are the five coordinate points and the line $\P(x_2,y)$. So the only possible singular points come from points in $Y\cap\P(x_2,y)$. This consists of three points: $(0:0:0:0:1),(0:1:0:0:1),(0:1:0:0:-1)$. For the two last points, we can use standard coordinates because these two points are not singular in $Y$ (nor in $\P(1,1,1,1,2)$). In fact, we can take ${x_1\over x_2},{x_3\over x_2},{x_4\over x_2}$ as affine coordinates. Regarding our action as the multiplication by $(\epsilon^4,1,\epsilon^2,\epsilon^3,1)$, it is clear that for the second and third points we get a quotient singularity ${1\over5}(4,2,3)={1\over 5}(1,2,3)$. Moreover, since $-K_X$ is represented locally by the divisor $\{{x_1\over x_2}=0\}$ and $\sigma$ by ${x_4^2\over x_2^2}$, we obtain that $\sigma$ is, locally, equal to the canonical divisor (i.e. $l_Q=1$ for both points). For the first point, we immediately see that we can take ${x_1\over y},{x_3\over y},{x_4\over y}$ as (analytically) local orbinates. Here, we have to consider the action in $\C^3$ which gives the singularity in $Y$: it takes $(a,b,c)$ to $(-a,-b,-c)$. Now we also have another action, generated by the morphism that takes $(a,b,c)$ to $(\epsilon^4a,\epsilon^2b,\epsilon^3c)$. This means that we actually have a $\Z/(10)$ action generated by the two morphisms just written (the opposite of a nontrivial 5th root of unit is a nontrivial 10th root of unit). So we get a singularity of type ${1\over10}(1\times5+4\times2,1\times5+2\times2,1\times5+3\times2)={1\over10}(3,9,1)={1\over10}(1,3,7)$. Working as before, we can see that $l_Q$ is 6 as expected. This means that the example we got does have only terminal singularities, and so it is a Fano--Enriques threefold.
\end{ex}

\vskip 10pt
\section{The search for Fano--Enriques threefolds}
In this section we list all non-special (in the sense of Remark \ref{nonspecial}) Fano--Enriques threefolds of codimension 1, 2 and 3. To this purpose we  give first some restrictions for the numerical data that gives (after a computational search), just 39 possibilities for $\Bt$. Then, we combine it with all possible covers (chosen from the lists of {\bf [F]} and {\bf [A]} or {\bf [GRDW]}) to give the quotients. The lists in this section have been found thanks to Magma (see {\bf [M]}) which was used also to do the computer search for Tables 1.$r$ ($r\in\{ 2,3,4,5,6,8\}$).

These are some immediate restrictions for the numerical data of a Fano--Enriques threefold:

\begin{enumerate}

\item from Bogomolov's instability, as said in {\bf [ABR]}, $-K_Xc_2>0$, so applying Orbifold-RR to the canonical class we get
$$\sum_{Q\in\B_X}\Big(r_Q-{1\over r_Q}\Big)<24.$$

\item the same for the cover $Y$:
$$\sum_{Q\in\B_Y}\Big(r_Q-{1\over r_Q}\Big)<24.$$

\item the Lefschetz number is $L(g,\O_Y)= \sum (-1)^i$Trace$(g^* |   _{H^i(\O_X)})$. By the Atiyah-Singer-Segal formula, it is a sum of various contributions over Fix$(g)$. If Fix$(g)$ = $\emptyset$, then $L(g,\O_Y)=0$. By Kodaira vanishing, since
the covering $Y$ is a Fano threefold, we get that $L(g,\O_Y)=$Trace$(g^* | _{H^0(\O_X)})$, and $g^* | _{H^0(\O_X)}$ is not zero. Therefore the action must have a fixed point. This means that, if $r$ is the order of the torsion divisor $\sigma$, we have a ${1\over kr}$ singularity for at least one $k\in \Z$ positive. So $r \le 24$ by (1).

\item if we call $\Bt$ the subset of $\B$ consisting on the singularities $Q$ where $\tau$ is not Cartier (i.e. $l_Q\ne 0$), it is clear that $\chi(-nK_X+\tau)-\chi(-nK_X)$ depends only on $\Bt$ and the coefficients $l_Q$, not on the rest of $\B$ and $-K_X^3$. These numbers must be integer for $\tau=i\sigma, i \in \{0,...,r-1\}$.  

\item Since $-K_X$ is ample and $\sigma$ is numerically trivial, $-K_X+i\sigma$ is ample for any $i\in\Z/(r)$. Therefore, by Kodaira's vanishing, $\chi(i\sigma)=h^0(i\sigma)\ge 0$. 
\end{enumerate}

After doing an exhaustive computer search (testing all baskets with all possible combinations of $l_Q$), we found that there are 39 possible $\Bt$ (with $r=2,3,4,5,6,8$) satisfying the restrictions. 
The notation we  use is
$$\Big({1\over r_Q}(1,a,-a)\Big)_{l_Q}.$$
 These are the results, which we divide according to the group $\Z/(r)$ acting

\centerline{\bf Table 1.2. Subsets $\Bt$ for order 2 actions:}

$$\Big({1\over 2}(1,1,1)\Big)_1,\Big({1\over 14}(1,1,13)\Big)_7 \eqno(\Bt 2.1)$$ 
$$\Big({1\over 2}(1,1,1)\Big)_1,\Big({1\over 14}(1,3,11)\Big)_7 \eqno(\Bt 2.2)$$
$$\Big({1\over 2}(1,1,1)\Big)_1,\Big({1\over 14}(1,5,9)\Big)_7 \eqno(\Bt 2.3)$$
$$\Big({1\over 4}(1,1,3)\Big)_2,\Big({1\over 12}(1,1,11)\Big)_6 \eqno(\Bt 2.4)$$
$$\Big({1\over 4}(1,1,3)\Big)_2,\Big({1\over 12}(1,5,7)\Big)_6 \eqno(\Bt 2.5)$$
$$\Big({1\over 6}(1,1,5)\Big)_3,\Big({1\over 10}(1,1,9)\Big)_5\eqno(\Bt 2.6)$$
$$\Big({1\over 6}(1,1,5)\Big)_3,\Big({1\over 10}(1,3,7)\Big)_5 \eqno(\Bt 2.7)$$
$$\Big({1\over 8}(1,1,7)\Big)_4,\Big({1\over 8}(1,1,7)\Big)_4 \eqno(\Bt 2.8)$$
$$\Big({1\over 8}(1,1,7)\Big)_4,\Big({1\over 8}(1,3,5)\Big)_4 \eqno(\Bt 2.9)$$
$$\Big({1\over 8}(1,3,5)\Big)_4,\Big({1\over 8}(1,3,5)\Big)_4 \eqno(\Bt 2.10)$$
$$3\times\Big({1\over 2}(1,1,1)\Big)_1,\Big({1\over 10}(1,1,9)\Big)_5 \eqno(\Bt 2.11)$$
$$3\times\Big({1\over 2}(1,1,1)\Big)_1,\Big({1\over 10}(1,3,7)\Big)_5 \eqno(\Bt 2.12)$$
$$2\times\Big({1\over 2}(1,1,1)\Big)_1,\Big({1\over 4}(1,1,3)\Big)_2,\Big({1\over 8}(1,1,7)\Big)_4 \eqno(\Bt 2.13)$$
$$2\times\Big({1\over 2}(1,1,1)\Big)_1,\Big({1\over 4}(1,1,3)\Big)_2,\Big({1\over 8}(1,3,5)\Big)_4 \eqno(\Bt 2.14)$$
$$2\times\Big({1\over 2}(1,1,1)\Big)_1,2\times\Big({1\over 6}(1,1,5)\Big)_3 \eqno(\Bt 2.15)$$
$$\Big({1\over 2}(1,1,1)\Big)_1,2\times\Big({1\over 4}(1,1,3)\Big)_2,\Big({1\over 6}(1,1,5)\Big)_3 \eqno(\Bt 2.16)$$
$$4\times\Big({1\over 4}(1,1,3)\Big)_2 \eqno(\Bt 2.17)$$
$$5\times\Big({1\over 2}(1,1,1)\Big)_1,\Big({1\over 6}(1,1,5)\Big)_3 \eqno(\Bt 2.18)$$
$$4\times\Big({1\over 2}(1,1,1)\Big)_1,2\times\Big({1\over 4}(1,1,3)\Big)_2 \eqno(\Bt 2.19)$$
$$8\times\Big({1\over 2}(1,1,1)\Big)_1 \eqno(\Bt 2.20)$$
\vskip 10pt

\centerline{\bf Table 1.3. Subsets $\Bt$ for order 3 actions:}

$$\Big({1\over9}(1,1,8)\Big)_3,\Big({1\over9}(1,1,8)\Big)_6 \eqno(\Bt 3.1)$$
$$\Big({1\over9}(1,2,7)\Big)_3,\Big({1\over9}(1,2,7)\Big)_6 \eqno(\Bt 3.2)$$
$$\Big({1\over9}(1,4,5)\Big)_3,\Big({1\over9}(1,4,5)\Big)_6 \eqno(\Bt 3.3)$$
$$2\times\Big({1\over3}(1,1,2)\Big)_1,\Big({1\over12}(1,5,7)\Big)_4 \eqno(\Bt 3.4)$$
$$\Big({1\over3}(1,1,2)\Big)_1,\Big({1\over3}(1,1,2)\Big)_2,\Big({1\over6}(1,1,5)\Big)_2,\Big({1\over6}(1,1,5)\Big)_4 \eqno(\Bt 3.5)$$
$$4\times\Big({1\over3}(1,1,2)\Big)_1,\Big({1\over6}(1,1,5)\Big)_4 \eqno(\Bt 3.6)$$
$$3\times\Big({1\over3}(1,1,2)\Big)_1,3\times\Big({1\over3}(1,1,2)\Big)_2 \eqno(\Bt 3.7)$$
\vskip 10pt

\centerline{\bf Table 1.4. Subsets $\Bt$ for order 4 actions:}

$$\Big({1\over4}(1,1,3)\Big)_2,2\times\Big({1\over8}(1,3,5)\Big)_2 \eqno(\Bt 4.1)$$ 
$$2\times\Big({1\over2}(1,1,1)\Big)_1,\Big({1\over4}(1,1,3)\Big)_1,\Big({1\over 12}(1,5,7)\Big)_9 \eqno(\Bt 4.2)$$
$$2\times\Big({1\over2}(1,1,1)\Big)_1,\Big({1\over8}(1,1,7)\Big)_2,\Big({1\over8}(1,1,7)\Big)_6 \eqno(\Bt 4.3)$$
$$2\times\Big({1\over2}(1,1,1)\Big)_1,\Big({1\over8}(1,3,5)\Big)_2,\Big({1\over8}(1,3,5)\Big)_6 \eqno(\Bt 4.4)$$
$$2\times\Big({1\over2}(1,1,1)\Big)_1,2\times\Big({1\over4}(1,1,3)\Big)_1,2\times\Big({1\over4}(1,1,3)\Big)_3 \eqno(\Bt 4.5)$$
\vskip 10pt

\centerline{\bf Table 1.5. Subsets $\Bt$ for order 5 actions:}

$$2\times\Big({1\over5}(1,2,3)\Big)_1,\Big({1\over10}(1,3,7)\Big)_6 \eqno(\Bt 5.1)$$
$$\Big({1\over5}(1,1,4)\Big)_1,\Big({1\over5}(1,1,4)\Big)_2,\Big({1\over5}(1,1,4)\Big)_3,\Big({1\over5}(1,1,4)\Big)_4 \eqno(\Bt 5.2)$$
$$\Big({1\over5}(1,1,4)\Big)_1,\Big({1\over5}(1,1,4)\Big)_4,\Big({1\over5}(1,2,3)\Big)_1,\Big({1\over5}(1,2,3)\Big)_4 \eqno(\Bt 5.3)$$
$$\Big({1\over5}(1,2,3)\Big)_1,\Big({1\over5}(1,2,3)\Big)_2,\Big({1\over5}(1,2,3)\Big)_3,\Big({1\over5}(1,2,3)\Big)_4 \eqno(\Bt 5.4)$$
\vskip 10pt

\centerline{\bf Table 1.6. Subsets $\Bt$ for order 6 actions:}

$$2\times\Big({1\over3}(1,1,2)\Big)_1,\Big({1\over4}(1,1,3)\Big)_2,\Big({1\over12}(1,5,7)\Big)_{10} \eqno(\Bt 6.1)$$
$$2\times\Big({1\over2}(1,1,1)\Big)_1,\Big({1\over3}(1,1,2)\Big)_1,\Big({1\over3}(1,1,2)\Big)_2,\Big({1\over6}(1,1,5)\Big)_1,\Big({1\over6}(1,1,5)\Big)_5 \eqno(\Bt 6.2)$$
\vskip 10pt

\centerline{\bf Table 1.8. Subsets $\Bt$ for order 8 actions:}

$$\Big({1\over2}(1,1,1)\Big)_1,\Big({1\over4}(1,1,3)\Big)_1,\Big({1\over8}(1,3,5)\Big)_3,\Big({1\over8}(1,3,5)\Big)_7 \eqno(\Bt 8.1)$$
\vskip 10pt

\begin{rem}
In the previous tables we omitted the redundant cases that are given by a multiple of the torsion divisor (i.e. same singularities but all $l_Q$ multiplied by an integer) because they represent the same $\Bt$ but the new torsion divisor is now a multiple of the former one. For instance, $2\times\Big({1\over5}(1,2,3)\Big)_3,\Big({1\over10}(1,3,7)\Big)_8$ is ($\Bt 5.1)$ considering $3\sigma$ instead of $\sigma$.
\end{rem}

Now we can use this list of possible subsets $\Bt$ to check which among all known Fano threefolds admit a Fano--Enriques quotient. We  check three lists, which we do not reproduce here due to their size. The first two, of codimension 1 and 2, due to Reid and Fletcher respectively, can be found in {\bf [F]}. The other one, of codimension 3, is due to Alt{\i}nok and is in {\bf [A]}.
For the first list (Reid's codimension 1) we get:

\begin{prop} 
Exactly 12 Fano--Enriques threefolds can be obtained as quotients of Reid's 95 Fano hypersurfaces. They have torsion 2,3 and 5 and are given in Table 2.
\end{prop}
\begin{pf} 
It is a case by case proof of the result. As a sample, we retake the case of Example \ref{fano} and Example \ref{ej1}: Let us suppose $\Z/(r)$ is acting on $Y_5\subset\P(1,1,1,1,2)$ (we explained in Remark \ref{obs-accion} that the action on $Y$ induces a diagonal action on $\P(1,1,1,1,2)$). Necessarily, (0:0:0:0:1) belongs to $Y$ and is a fixed point by the action on $\P(1,1,1,1,2)$, since it is the only singular point on this weighted projective space. Therefore, a singularity of type ${1\over 2r}$  appears in the quotient $X$. In fact, the basket has to be $\B =\Bt=\{ {1\over r_1},...,{1\over r_l}, {1\over 2r}\}$, with $r_i|r$. No ($\Bt r.i$) satisfies this condition for $r=2,4,6,8$. There is a possibility for $\Z/(3)$, with the basket $\{4\times{1\over3},{1\over6}\}$ ($\Bt$ 3.6), but it is also impossible because it should be $-K_X^3={5\over 2r}={5\over6}$, which is in contradiction with Remark \ref{selfint}, which implies that, for this basket, $-K_X^3={3\over2}+2k, \ k\in\Z$. Therefore, $r=5$ and the only possible quotient is the example we already know (it is No. 2 in Table 2 below). Now we would play the game we played in Example \ref{ej1} to observe that this is the only possibility with these invariants.
\end{pf}

\begin{rem}
For the following tables, we  list these data for each element:
\begin{itemize}
\item the cover, which  is a complete intersection $Y_{d_1,...,d_l}\subset \P(a_1,...,a_n)$ of hypersurfaces of degrees $d_1,...,d_l$
\item the action on the cover
\item $\Bt:=\{$singularities in the basket where the torsion divisor is not trivial with their respective coefficients $l_Q$ $\}$.
\item $\B \backslash \Bt:=$ rest of the basket
\end{itemize}
\end{rem}

\vskip 10pt
\centerline{\bf Table 2. Fano--Enriques threefolds from codimension 1 Fano threefolds:}
\vskip 10pt
\leftline{\bf No. 1}
\begin{itemize}
\item {\bf cover:} $Y_{4}\subset \P(1,1,1,1,1)$
\item {\bf action:} $\Z/(5)$ acts by $(1,\epsilon,\epsilon^2,\epsilon^3,\epsilon^4)$, $\epsilon=e^{{2\pi\over 5}}$
\item $\Bt=\{\Big({1\over 5}(1,2,3)\Big)_1,\Big({1\over 5}(1,2,3)\Big)_2,\Big({1\over 5}(1,2,3)\Big)_3,\Big({1\over 5}(1,2,3)\Big)_4\}=(\Bt 5.4)$
\item $\B \backslash \Bt=\emptyset$
\end{itemize}
\leftline{\bf No. 2}
\begin{itemize}
\item {\bf cover:} $Y_{5}\subset \P(1,1,1,1,2)$
\item {\bf action:} $\Z/(5)$ acts by $(1,\epsilon,\epsilon^3,\epsilon^4,\epsilon^2)$, $\epsilon=e^{{2\pi\over 5}}$
\item $\Bt=\{2\times \Big({1\over 5}(1,2,3)\Big)_1, \Big({1\over 10}(1,3,7)\Big)_6\}=(\Bt 5.1)$
\item $\B \backslash \Bt=\emptyset$
\end{itemize}
\leftline{\bf No. 3}
\begin{itemize}
\item {\bf cover:} $Y_{6}\subset\P(1,1,1,2,2)$
\item {\bf action:} $\Z/(3)$ acts by $(1,\epsilon,\epsilon^2,\epsilon,\epsilon^2)$, $\epsilon=e^{{2\pi\over 3}}$
\item $\Bt=\{3\times\Big({1\over3}(1,1,2)\Big)_1,3\times\Big({1\over3}(1,1,2)\Big)_2\}=(\Bt 3.7)$ 
\item $\B \backslash \Bt=\{ {1\over 2}(1,1,1)\}$
\end{itemize}
\leftline{\bf No. 4}
\begin{itemize}
\item {\bf cover:} $Y_{8}\subset\P(1,1,1,2,4)$
\item {\bf action:} $\Z/(2)$ acts by $(+,-,-,-,-)$
\item $\Bt=\{8\times\Big({1\over 2}(1,1,1)\Big)_1\}=(\Bt 2.20)$ 
\item $\B \backslash \Bt=\{ {1\over 2}(1,1,1)\}$
\end{itemize}
\leftline{\bf No. 5}
\begin{itemize}
\item {\bf cover:} $Y_{9}\subset\P(1,1,1,3,4)$
\item {\bf action:} $\Z/(3)$ acts by $(1,\epsilon,\epsilon^2,\epsilon^2,\epsilon)$, $\epsilon=e^{{2\pi\over 3}}$
\item $\Bt=\{2\times\Big({1\over 3}(1,1,2)\Big)_1,\Big({1\over 12}(1,5,7)\Big)_4\}=(\Bt 3.4)$ 
\item $\B \backslash \Bt=\emptyset$
\end{itemize}
\leftline{\bf No. 6}
\begin{itemize}
\item {\bf cover:} $Y_{9}\subset\P(1,1,2,3,3)$
\item {\bf action:} $\Z/(3)$ acts by $(1,\epsilon,\epsilon^2,\epsilon,\epsilon^2)$, $\epsilon=e^{{2\pi\over 3}}$
\item $\Bt=\{4\times\Big({1\over 3}(1,1,2)\Big)_1,\Big({1\over 6}(1,1,5)\Big)_4\}=(\Bt 3.6)$ 
\item $\B \backslash \Bt=\{ {1\over 3}(1,1,2)\}$
\end{itemize}
\leftline{\bf No. 7}
\begin{itemize}
\item {\bf cover:} $Y_{12)}\subset\P(1,1,2,3,6)$
\item {\bf action:} $\Z/(2)$ acts by $(+,-,-,-,-)$
\item $\Bt=\{4\times\Big({1\over2}(1,1,1)\Big)_1,2\times\Big({1\over4}(1,1,3)\Big)_2\}=(\Bt 2.19)$ 
\item $\B \backslash \Bt=\{ {1\over 3}(1,1,2)\}$
\end{itemize}
\leftline{\bf No. 8}
\begin{itemize}
\item {\bf cover:} $Y_{14}\subset\P(1,1,2,4,7)$
\item {\bf action:} $\Z/(2)$ acts by $(+,-,-,-,-)$
\item $\Bt=\{2\times\Big({1\over 2}(1,1,1)\Big)_1,\Big({1\over 4}(1,1,3)\Big)_2,\Big({1\over 8}(1,3,5)\Big)_4\}=(\Bt 2.14)$ 
\item $\B \backslash \Bt=\{ {1\over 2}(1,1,1)\}$
\end{itemize}
\leftline{\bf No. 9}
\begin{itemize}
\item {\bf cover:} $Y_{16}\subset\P(1,1,2,5,8)$
\item {\bf action:} $\Z/(2)$ acts by $(+,-,-,-,-)$
\item $\Bt=\{3\times\Big({1\over 2}(1,1,1)\Big)_1,\Big({1\over 10}(1,3,)\Big)_5\}=(\Bt 2.12)$ 
\item $\B \backslash \Bt=\{ {1\over 2}(1,1,1)\}$
\end{itemize}
\leftline{\bf No. 10}
\begin{itemize}
\item {\bf cover:} $Y_{16}\subset\P(1,1,3,4,8)$
\item {\bf action:} $\Z/(2)$ acts by $(+,-,-,-,-)$
\item $\Bt=\{5\times\Big({1\over 2}(1,1,1)\Big)_1,\Big({1\over 6}(1,1,5)\Big)_3\}=(\Bt 2.18)$ 
\item $\B \backslash \Bt=\{ {1\over 4}(1,1,3)\}$
\end{itemize}
\leftline{\bf No. 11}
\begin{itemize}
\item {\bf cover:} $Y_{20}\subset\P(1,2,3,5,10)$
\item {\bf action:} $\Z/(2)$ acts by $(+,-,-,-,-)$
\item $\Bt=\{\Big({1\over 2}(1,1,1)\Big)_1,2\times\Big({1\over 4}(1,1,3)\Big)_2,\Big({1\over 6}(1,1,5)\Big)_3\}=(\Bt 2.16)$ 
\item $\B \backslash \Bt=\{ {1\over 5}(1,2,3)\}$
\end{itemize}
\leftline{\bf No. 12}
\begin{itemize}
\item {\bf cover:} $Y_{24}\subset\P(1,2,3,7,12)$
\item {\bf action:} $\Z/(2)$ acts by $(+,-,-,-,-)$
\item $\Bt=\{\Big({1\over 2}(1,1,1)\Big)_1,\Big({1\over 14}(1,5,9)\Big)_7\}=(\Bt 2.3)$ 
\item $\B \backslash \Bt=\{{1\over 2}(1,1,1) ,{1\over 3}(1,1,2)\}$
\end{itemize}

\begin{rem}\label{2nd-degree}
Equations defining $Y_d$ in the list have to be in an $H^0(-nK_Y+i\sigma)$. The degree $n$ is given by the corresponding $d$ subindex, so we have to complete the data giving $i$. All cases above have second degree $0\in\Z/(r)$.
\end{rem}

\begin{rem}\label{algor}
Although we have shown how to get by hand all Fano--Enriques quotients from a list of Fano threefolds, there is a general way to get them using a computer. The algorithm is divided into the following steps, which we should repeat for every ($\Bt\_.\_$) and every candidate $Y$ to be a cover:
\begin{enumerate}
\item for each singularity $Q\in \Bt$, of type $\Big({1\over r_Q}(1,a_Q,-a_Q)\Big)_{l_Q}$ we describe its preimage by the covering, which  consists in a precise number of points with the same determined type of singularity. In fact, if $d_Q={\mathrm gcd}(r_Q,l_Q)$, then the order of the torsion divisor $\sigma$ in the singularty $Q$ is $\alpha_Q={r_Q\over d_Q}$. This means that $Q$ is covered by ${r\over \alpha_Q}$ points where the local Picard group is $\Z/({r_Q\over d_Q})$ (of course, this is a ${1\over d_Q}$ singularity). So we have to find out the numbers $d_Q,a,b,c$ that describe the ${1\over d_Q}(a,b,c)$ quotient singularity of the points in the preimage of $Q$. We can simplify it considering a ${1\over d_Q}(a,b,c)$ singularity under a $\Z/(\beta)$ action, where $\beta$ is a multiple of $\alpha_Q$ (see Remark \ref{obs1}). We can suppose that the local orbinates have an action of type ${1\over \beta}(x,y,z)$. Since $r_Q=\beta d_Q$, composing the automorphisms generating the representations of $\Z/(d_Q)$, we get that $Q$ is a ${1\over \beta d_Q}(xd_Q+a\beta, yd_Q+b\beta, zd_Q+c\beta)$ singularity (even if $\beta$ and $d_Q$ are not coprime, in which case, all numbers can be reduced) and then, $a,b$ and $c$ are completely determined by $1,a_Q$ and $-a_Q$ (because $Q$ is a singularity of type ${1\over r_Q}(1,a_Q,-a_Q)$). But looking at all 39 subsets $\Bt$ that are listed before, the singularity in the cover $Y$ has to be terminal. Checking all possible subsets $\Bt$ it is obvious that $d_Q$ has to be 2, 3, 4, 5, 6 or 7 (or 1, but then it is a regular point and we do not consider this case since it does not affect the numerical data). There is only one quotient singularity of order $d_Q \in \{2,3,4,6\}$ which is ${1\over d_Q}(1,1,-1)$. So we need to study the cases $d_Q=5,7$. But then $\alpha_Q=\beta=2$, so $x=y=z=1$ (in other case, $a_Q$ would not be coprime with $r_Q$) and the only thing remaining is to solve the equation ${1\over \beta d_Q}(xd_Q+a\beta, yd_Q+b\beta, zd_Q+c\beta)={1\over r_Q}(1,a_Q,-a_Q)$.
\item repeating this for every singularity in $\Bt$, we have a completely determined subset $\tilde\B$ of the basket $\B_Y$ consisting of all singularities where the action is not free. First condition for $Y$ to be a cover of a Fano--Enriques $X$ with the selected $\Bt$ is that $\tilde\B$ is contained in $\B_Y$. 
\item for the rest of $\B_Y$, the action is free on every member, so the second condition is that $\B_Y\backslash \tilde\B$ is divided in orbits, each of them consisting of $r$ singularities of the same type. Chosing an element of each orbit we  get the set $\B_X\backslash\Bt$ of all singularities in $\B_X$ where the torsion divisor is trivial.
\item from Remark \ref{selfint}, we know that, for a Fano threefold, 
$$-K_X^3=2k+\sum_{Q\in\B}{b_Q(r-b_Q)\over r_Q}.$$

It is obvious that $K_X^3={1\over r}K_Y^3$, so the last condition to test for the numerical data is that
$${K_Y^3\over r}-\sum_{Q\in\B_X}{b_Q(r-b_Q)\over r_Q}=2k.$$
\item for a pair $Y,\Bt$ satisfying the conditions in (2), (3) and (4), the numerical data for a Fano--Enriques $X$ covered by $Y$ associated to $\Bt$  then consists of $-K_X^3=-{1\over r}K_Y^3$ and a basket $\B$ that is the union of $\Bt$ and $\B_X\backslash\Bt$ constructed as in (3)
\item Finally, we apply the graded rings method (see Example \ref{ej1}) to the numerical data obtained in (5). Observe that even if we start from a list of Fano threefolds of a given codimension, this final step could produce a Fano--Enriques $X$ with an unexpected cover  of higher codimension (see Remark \ref{ej2}). We  divide our lists by means of the codimension of the resulting cover, not in terms of the codimension of the starting Fano threefold.
\end{enumerate}
\end{rem}

\begin{rem}\label{obs1}
In step (1) of the algorithm in Remark \ref{algor} we have to take $\beta$ to be a multiple of $\alpha$ even if they are actually equal. This is because, in order to use the expression ${1\over \beta d_Q}(xd_Q+a\beta, yd_Q+b\beta, zd_Q+c\beta)$,  we  need to use a fake $\beta$.  Let us think for instance of a $\Z/(2)$ action that fixes a singularity of type ${1\over 2}(1,1,1)$. For example we can muliply by $(+,+,-,-)$ in $\P(1,1,1,2)$. To dehomogeneize by the last coordinate, we have to think of the action as multiplication by $(i,i,-i,1)$, so $\beta$ seems to be 4 although it is actually 2. We  work as if $\beta=4$ and get a singularity of type ${1\over 8}(1\cdot 2+1\cdot \beta,1\cdot 2+1\cdot \beta,3\cdot 2+1\cdot \beta)$ (i.e. a singularity of type ${1\over 4}(1,1,3)$).
\end{rem}

\vskip 10pt

Now we use the algorithm explained in Remark \ref{algor} to test the lists of Fano threefolds of codimension one, two and three, obtained respectively by Reid, Fletcher and Alt\i nok to find candidates to be a cover of a Fano--Enriques threefold. Of course, for codimension one we reobtain Table 2. We list in the next tables the cases of codimension two and three, keeping the notations of Table 2.

\vskip 10pt
\centerline{\bf Table 3: Fano--Enriques threefolds from codimension 2 Fano threefolds:}
\vskip 10pt
\leftline{\bf No. 1a}
\begin{itemize}
\item {\bf cover:} $Y_{2,3}\subset\P(1,1,1,1,1,1)$
\item {\bf action:} $\Z/(3)$ acts by $(1,1,\epsilon,\epsilon,\epsilon^2,\epsilon^2)$
\item $\Bt=\{3\times\Big({1\over 3}(1,1,2)\Big)_1,3\times\Big({1\over 3}(1,1,2)\Big)_2\}=(\Bt 3.7)$ 
\item $\B \backslash \Bt=\emptyset$
\end{itemize}
\leftline{\bf No. 1b}
\begin{itemize}
\item {\bf cover:} $Y_{2,3}\subset\P(1,1,1,1,1,1)$
\item {\bf action:} $\Z/(5)$ acts by $(1,1,\epsilon,\epsilon^2,\epsilon^3,\epsilon^4)$
\item $\Bt=\{\Big({1\over 5}(1,1,4)\Big)_1,\Big({1\over 5}(1,1,4)\Big)_2,\Big({1\over 5}(1,1,4)\Big)_3,\Big({1\over 5}(1,1,4)\Big)_4\}=(\Bt 5.2)$ 
\item $\B \backslash \Bt=\emptyset$
\end{itemize}
\leftline{\bf No. 2}
\begin{itemize}
\item {\bf cover:} $Y_{3,3}\subset\P(1,1,1,1,1,2)$
\item {\bf action:} $\Z/(3)$ acts by $(1,1,\epsilon,\epsilon,\epsilon^2,\epsilon^2)$
\item $\Bt=\{4\times\Big({1\over 3}(1,1,2)\Big)_1,\Big({1\over 6}(1,1,5)\Big)_4\}=(\Bt 3.6)$ 
\item $\B \backslash \Bt=\emptyset$
\end{itemize}
\leftline{\bf No. 3}
\begin{itemize}
\item {\bf cover:} $Y_{2,4}\subset\P(1,1,1,1,1,2)$
\item {\bf action:} $\Z/(2)$ acts by $(+,+,-,-,-,-)$
\item $\Bt=\{8\times\Big({1\over 2}(1,1,1)\Big)_1\}=(\Bt 2.20)$ 
\item $\B \backslash \Bt=\emptyset$
\end{itemize}
\leftline{\bf No. 4a}
\begin{itemize}
\item {\bf cover:} $Y_{3,4}\subset\P(1,1,1,1,2,2)$
\item {\bf action:} $\Z/(2)$ acts by $(+,+,-,-,-,-)$
\item $\Bt=\{4\times\Big({1\over 2}(1,1,1)\Big)_1,2\times\Big({1\over 4}(1,1,3)\Big)_2\}=(\Bt 2.19)$ 
\item $\B \backslash \Bt=\emptyset$
\end{itemize}
\leftline{\bf No. 4b}
\begin{itemize}
\item {\bf cover:} $Y_{3,4}\subset\P(1,1,1,1,2,2)$
\item {\bf action:} $\Z/(3)$ acts by $(1,1,\epsilon,\epsilon^2,\epsilon,\epsilon^2)$
\item $\Bt=\{\Big({1\over 3}(1,1,2)\Big)_1,\Big({1\over 3}(1,1,2)\Big)_2,\Big({1\over 6}(1,1,5)\Big)_2,\Big({1\over 6}(1,1,5)\Big)_4\}=(\Bt 3.5)$ 
\item $\B \backslash \Bt=\emptyset$
\end{itemize}
\leftline{\bf No. 4c}
\begin{itemize}
\item {\bf cover:} $Y_{3,4}\subset\P(1,1,1,1,2,2)$
\item {\bf action:} $\Z/(4)$ acts by $(1,1,\epsilon,\epsilon^3,\epsilon,\epsilon^3)=(+,+,i,-i,i,-i)$
\item $\Bt=\{2\times\Big({1\over 2}(1,1,1)\Big)_1,\Big({1\over 8}(1,1,7)\Big)_2,\Big({1\over 8}(1,1,7)\Big)_6\}=(\Bt 4.4)$ 
\item $\B \backslash \Bt=\emptyset$
\end{itemize}
\leftline{\bf No. 4d}
\begin{itemize}
\item {\bf cover:} $Y_{3,4}\subset\P(1,1,1,1,2,2)$
\item {\bf action:} $\Z/(4)$ acts by $(1,\epsilon,\epsilon^2,\epsilon^3,\epsilon,\epsilon^3)=(+,i,-,-i,i,-i)$
\item $\Bt=\{2\times\Big({1\over 2}(1,1,1)\Big)_1,\Big({1\over 8}(1,3,5)\Big)_2,\Big({1\over 8}(1,3,5)\Big)_6\}=(\Bt 4.4)$ 
\item $\B \backslash \Bt=\emptyset$
\end{itemize}
\leftline{\bf No. 5a}
\begin{itemize}
\item {\bf cover:} $Y_{4,4}\subset\P(1,1,1,1,2,3)$
\item {\bf action:} $\Z/(4)$ acts by $(+,+,-,-,-,-)$
\item $\Bt=\{5\times\Big({1\over 2}(1,1,1)\Big)_1,\Big({1\over 6}(1,1,5)\Big)_3\}=(\Bt 2.18)$ 
\item $\B \backslash \Bt=\emptyset$
\end{itemize}
\leftline{\bf No. 5b}
\begin{itemize}
\item {\bf cover:} $Y_{4,4}\subset\P(1,1,1,1,2,3)$
\item {\bf action:} $\Z/(4)$ acts by $(+,i,-,-i,i,-i)$
\item $\Bt=\{2\times\Big({1\over 2}(1,1,1)\Big)_1,\Big({1\over 4}(1,1,3)\Big)_1,\Big({1\over 12}(1,5,7)\Big)_9\}=(\Bt 4.2)$ 
\item $\B \backslash \Bt=\emptyset$
\end{itemize}
\leftline{\bf No. 6a}
\begin{itemize}
\item {\bf cover:} $Y_{4,4}\subset\P(1,1,1,2,2,2)$
\item {\bf action:} $\Z/(2)$ acts by $(+,-,-,+,-,-)$
\item $\Bt=\{8\times\Big({1\over 2}(1,1,1)\Big)_1\}=(\Bt 2.20)$ 
\item $\B \backslash \Bt=\{2\times{1\over 2}(1,1,1)\}$
\end{itemize}
\leftline{\bf No. 6b}
\begin{itemize}
\item {\bf cover:} $Y_{4,4}\subset\P(1,1,1,2,2,2)$
\item {\bf action:} $\Z/(2)$ acts by $(+,+,-,-,-,-)$
\item $\Bt=\{4\times\Big({1\over 4}(1,1,3)\Big)_2\}=(\Bt 2.17)$ 
\item $\B \backslash \Bt=\emptyset$
\end{itemize}
\leftline{\bf No. 6c}
\begin{itemize}
\item {\bf cover:} $Y_{4,4}\subset\P(1,1,1,2,2,2)$
\item {\bf action:} $\Z/(4)$ acts by $(+,i,-i,i,-,-i)$
\item $\Bt=\{2\times\Big({1\over 2}(1,1,1)\Big)_1,2\times\Big({1\over 4}(1,1,3)\Big)_1,4\times\Big({1\over 4}(1,1,3)\Big)_3\}=(\Bt 4.5)$ 
\item $\B \backslash \Bt=\{{1\over 2}(1,1,1)\}$
\end{itemize}
\leftline{\bf No. 6d}
\begin{itemize}
\item {\bf cover:} $Y_{4,4}\subset\P(1,1,1,2,2,2)$
\item {\bf action:} $\Z/(4)$ acts by $(+,-,-i,i,i,-i)$
\item $\Bt=\{\Big({1\over 4}(1,1,3)\Big)_2,2\times\Big({1\over 8}(1,3,5)\Big)_2\}=(\Bt 4.1)$ 
\item $\B \backslash \Bt=\emptyset$
\end{itemize}
\leftline{\bf No. 7}
\begin{itemize}
\item {\bf cover:} $Y_{4,5}\subset\P(1,1,1,2,2,3)$
\item {\bf action:} $\Z/(2)$ acts by $(+,+,-,-,-,-)$
\item $\Bt=\{\Big({1\over 2}(1,1,1)\Big)_1,2\times\Big({1\over 4}(1,1,3)\Big)_2,\Big({1\over 6}(1,1,5)\Big)_3\}=(\Bt 2.16)$ 
\item $\B \backslash \Bt=\emptyset$
\end{itemize}
\leftline{\bf No. 8a}
\begin{itemize}
\item {\bf cover:} $Y_{4,6}\subset\P(1,1,1,2,3,3)$
\item {\bf action:} $\Z/(2)$ acts by $(+,-,-,-,+,-)$
\item $\Bt=\{8\times\Big({1\over 2}(1,1,1)\Big)_1\}=(\Bt 2.20)$ 
\item $\B \backslash \Bt=\{{1\over 2}(1,1,1)\}$
\end{itemize}
\leftline{\bf No. 8b}
\begin{itemize}
\item {\bf cover:} $Y_{4,6}\subset\P(1,1,1,2,3,3)$
\item {\bf action:} $\Z/(2)$ acts by $(+,+,-,-,-,-)$
\item $\Bt=\{2\times\Big({1\over 2}(1,1,1)\Big)_1,2\times\Big({1\over 6}(1,1,5)\Big)_3\}=(\Bt 2.15)$ 
\item $\B \backslash \Bt=\emptyset$
\end{itemize}
\leftline{\bf No. 8c}
\begin{itemize}
\item {\bf cover:} $Y_{4,6}\subset\P(1,1,1,2,3,3)$
\item {\bf action:} $\Z/(3)$ acts by $(1,\epsilon,\epsilon^2,1,\epsilon,\epsilon^2)$
\item $\Bt=\{\Big({1\over 9}(1,2,7)\Big)_3,\Big({1\over 9}(1,2,7)\Big)_6\}=(\Bt 3.2)$ 
\item $\B \backslash \Bt=\emptyset$
\end{itemize}
\leftline{\bf No. 9}
\begin{itemize}
\item {\bf cover:} $Y_{5,6}\subset\P(1,1,1,2,3,4)$
\item {\bf action:} $\Z/(2)$ acts by $(+,+,-,-,-,-)$
\item $\Bt=\{2\times\Big({1\over 2}(1,1,1)\Big)_1,\Big({1\over 4}(1,1,3)\Big)_2,\Big({1\over 8}(1,1,7)\Big)_4\}=(\Bt 2.13)$ 
\item $\B \backslash \Bt=\emptyset$
\end{itemize}
\leftline{\bf No. 10}
\begin{itemize}
\item {\bf cover:} $Y_{6,8}\subset\P(1,1,1,3,4,5)$
\item {\bf action:} $\Z/(2)$ acts by $(+,+,-,-,-,-)$
\item $\Bt=\{3\times\Big({1\over 2}(1,1,1)\Big)_1,\Big({1\over 10}(1,1,9)\Big)_5\}=(\Bt 2.11)$ 
\item $\B \backslash \Bt=\emptyset$
\end{itemize}
\leftline{\bf No. 11}
\begin{itemize}
\item {\bf cover:} $Y_{4,6}\subset\P(1,1,2,2,2,3)$
\item {\bf action:} $\Z/(2)$ acts by $(+,-,+,-,-,-)$
\item $\Bt=\{4\times\Big({1\over 2}(1,1,1)\Big)_1,2\times\Big({1\over 4}(1,1,3)\Big)_2\}=(\Bt 2.19)$ 
\item $\B \backslash \Bt=\{2\times{1\over 2}(1,1,1)\}$
\end{itemize}
\leftline{\bf No. 12}
\begin{itemize}
\item {\bf cover:} $Y_{5,6}\subset\P(1,1,2,2,3,3)$
\item {\bf action:} $\Z/(3)$ acts by $(1,1,\epsilon,\epsilon^2,\epsilon,\epsilon^2)$
\item $\Bt=\{\Big({1\over 9}(1,1,8)\Big)_3,\Big({1\over 9}(1,1,8)\Big)_6\}=(\Bt 3.1)$ 
\item $\B \backslash \Bt=\{{1\over 2}(1,1,1)\}$
\end{itemize}
\leftline{\bf No. 13}
\begin{itemize}
\item {\bf cover:} $Y_{4,8}\subset\P(1,1,2,2,3,4)$
\item {\bf action:} $\Z/(2)$ acts by $(+,-,+,-,-,-)$
\item $\Bt=\{5\times\Big({1\over 2}(1,1,1),\Big)_1,\Big({1\over 6}(1,1,5)\Big)_3\}=(\Bt 2.18)$ 
\item $\B \backslash \Bt=\{2\times{1\over 2}(1,1,1)\}$
\end{itemize}
\leftline{\bf No. 14}
\begin{itemize}
\item {\bf cover:} $Y_{6,7}\subset\P(1,1,2,3,3,4)$
\item {\bf action:} $\Z/(2)$ acts by $(+,-,-,+,-,-)$
\item $\Bt=\{2\times\Big({1\over 2}(1,1,1),\Big)_1,\Big({1\over 4}(1,1,3)\Big)_2,\Big({1\over 8}(1,3,5)\Big)_4\}=(\Bt 2.14)$ 
\item $\B \backslash \Bt=\{{1\over 3}(1,1,2)\}$
\end{itemize}
\leftline{\bf No. 15}
\begin{itemize}
\item {\bf cover:} $Y_{6,8}\subset\P(1,1,2,3,3,5)$
\item {\bf action:} $\Z/(2)$ acts by $(+,-,-,+,-,-)$
\item $\Bt=\{3\times\Big({1\over 2}(1,1,1)\Big)_1,\Big({1\over 10}(1,3,7)\Big)_5\}=(\Bt 2.11)$ 
\item $\B \backslash \Bt=\{{1\over 3}(1,1,2)\}$
\end{itemize}
\leftline{\bf No. 16a}
\begin{itemize}
\item {\bf cover:} $Y_{6,8}\subset\P(1,1,2,3,4,4)$
\item {\bf action:} $\Z/(2)$ acts by $(+,-,-,-,+,-)$
\item $\Bt=\{4\times\Big({1\over 2}(1,1,1)\Big)_1,2\times\Big({1\over 4}(1,1,3)\Big)_2\}=(\Bt 2.19)$ 
\item $\B \backslash \Bt=\{{1\over 4}(1,1,3)\}$
\end{itemize}
\leftline{\bf No. 16b}
\begin{itemize}
\item {\bf cover:} $Y_{6,8}\subset\P(1,1,2,3,4,4)$
\item {\bf action:} $\Z/(2)$ acts by $(+,+,-,-,-,-)$
\item $\Bt=\{2\times\Big({1\over 8}(1,1,7)\Big)_4\}=(\Bt 2.8)$ 
\item $\B \backslash \Bt=\{{1\over 2}(1,1,1)\}$
\end{itemize}
\leftline{\bf No. 16c}
\begin{itemize}
\item {\bf cover:} $Y_{6,8}\subset\P(1,1,2,3,4,4)$
\item {\bf action:} $\Z/(2)$ acts by $(+,-,-,+,-,-)$
\item $\Bt=\{2\times\Big({1\over 8}(1,3,5)\Big)_4\}=(\Bt 2.10)$ 
\item $\B \backslash \Bt=\{{1\over 2}(1,1,1)\}$
\end{itemize}
\leftline{\bf No. 17}
\begin{itemize}
\item {\bf cover:} $Y_{7,8}\subset\P(1,1,2,3,4,5)$
\item {\bf action:} $\Z/(2)$ acts by $(+,+,-,-,-,-)$
\item $\Bt=\{\Big({1\over 6}(1,1,5)\Big)_3,\Big({1\over 10}(1,1,9)\Big)_5\}=(\Bt 2.6)$ 
\item $\B \backslash \Bt=\{{1\over 2}(1,1,1)\}$
\end{itemize}
\leftline{\bf No. 18}
\begin{itemize}
\item {\bf cover:} $Y_{8,9}\subset\P(1,1,2,3,4,7)$
\item {\bf action:} $\Z/(2)$ acts by $(+,+,-,-,-,-)$
\item $\Bt=\{\Big({1\over 2}(1,1,1)\Big)_1,\Big({1\over 14}(1,3,11)\Big)_7\}=(\Bt 2.2)$ 
\item $\B \backslash \Bt=\{{1\over 2}(1,1,1)\}$
\end{itemize}
\leftline{\bf No. 19a}
\begin{itemize}
\item {\bf cover:} $Y_{8,10}\subset\P(1,1,2,4,5,6)$
\item {\bf action:} $\Z/(2)$ acts by $(+,+,-,-,-,-)$
\item $\Bt=\{\Big({1\over 4}(1,1,3)\Big)_2,\Big({1\over 12}(1,1,11)\Big)_6\}=(\Bt 2.4)$ 
\item $\B \backslash \Bt=\{{1\over 2}(1,1,1)\}$
\end{itemize}
\leftline{\bf No. 19b}
\begin{itemize}
\item {\bf cover:} $Y_{8,10}\subset\P(1,1,2,4,5,6)$
\item {\bf action:} $\Z/(2)$ acts by $(+,-,-,-,+,-)$
\item $\Bt=\{\Big({1\over 4}(1,1,3)\Big)_2,\Big({1\over 12}(1,5,7)\Big)_6\}=(\Bt 2.5)$ 
\item $\B \backslash \Bt=\{{1\over 2}(1,1,1)\}$
\end{itemize}
\leftline{\bf No. 20}
\begin{itemize}
\item {\bf cover:} $Y_{8,10}\subset\P(1,1,3,4,5,5)$
\item {\bf action:} $\Z/(2)$ acts by $(+,-,-,-,+,-)$
\item $\Bt=\{5\times\Big({1\over 2}(1,1,1)\Big)_1,\Big({1\over 6}(1,1,5)\Big)_3\}=(\Bt 2.18)$ 
\item $\B \backslash \Bt=\{{1\over 5}(1,1,5)\}$
\end{itemize}
\leftline{\bf No. 21}
\begin{itemize}
\item {\bf cover:} $Y_{10,12}\subset\P(1,1,3,5,6,7)$
\item {\bf action:} $\Z/(2)$ acts by $(+,+,-,-,-,-)$
\item $\Bt=\{\Big({1\over 2}(1,1,1)\Big)_1,\Big({1\over 14}(1,1,13)\Big)_7\}=(\Bt 2.1)$ 
\item $\B \backslash \Bt=\{{1\over 3}(1,1,2)\}$
\end{itemize}
\leftline{\bf No. 22a}
\begin{itemize}
\item {\bf cover:} $Y_{6,8}\subset\P(1,2,2,3,3,4)$
\item {\bf action:} $\Z/(2)$ acts by $(-,+,-,+,-,-)$
\item $\Bt=\{8\times\Big({1\over 2}(1,1,1)\Big)_1\}=(\Bt 2.20)$ 
\item $\B \backslash \Bt=\{3\times{1\over 2}(1,1,1),{1\over 3}(1,1,2)\}$
\end{itemize}
\leftline{\bf No. 22b}
\begin{itemize}
\item {\bf cover:} $Y_{6,8}\subset\P(1,2,2,3,3,4)$
\item {\bf action:} $\Z/(2)$ acts by $(+,+,-,-,-,-)$
\item $\Bt=\{2\times\Big({1\over 2}(1,1,1)\Big)_1,2\times\Big({1\over 6}(1,1,5)\Big)_3\}=(\Bt 2.15)$ 
\item $\B \backslash \Bt=\{3\times{1\over 2}(1,1,1)\}$
\end{itemize}
\leftline{\bf No. 22c}
\begin{itemize}
\item {\bf cover:} $Y_{6,8}\subset\P(1,2,2,3,3,4)$
\item {\bf action:} $\Z/(2)$ acts by $(+,-,-,+,-,-)$
\item $\Bt=\{4\times\Big({1\over 4}(1,1,3)\Big)_2\}=(\Bt 2.17)$ 
\item $\B \backslash \Bt=\{{1\over 2}(1,1,1),{1\over 3}(1,1,2)\}$
\end{itemize}
\leftline{\bf No. 22d}
\begin{itemize}
\item {\bf cover:} $Y_{6,8}\subset\P(1,2,2,3,3,4)$
\item {\bf action:} $\Z/(3)$ acts by $(1,\epsilon,\epsilon^2,\epsilon,\epsilon^2,1)$
\item $\Bt=\{\Big({1\over 9}(1,4,5)\Big)_3,\Big({1\over 9}(1,4,5)\Big)_6\}=(\Bt 3.3)$ 
\item $\B \backslash \Bt=\{2\times{1\over 2}(1,1,1)\}$
\end{itemize}
\leftline{\bf No. 23}
\begin{itemize}
\item {\bf cover:} $Y_{6,10}\subset\P(1,2,2,3,4,5)$
\item {\bf action:} $\Z/(2)$ acts by $(+,+,-,-,-,-)$
\item $\Bt=\{2\times\Big({1\over 2}(1,1,1)\Big)_1,\Big({1\over 4}(1,1,3)\Big)_2,\Big({1\over 8}(1,1,7)\Big)_4\}=(\Bt 2.13)$ 
\item $\B \backslash \Bt=\{3\times{1\over 2}(1,1,1)\}$
\end{itemize}
\leftline{\bf No. 24}
\begin{itemize}
\item {\bf cover:} $Y_{8,9}\subset\P(1,2,3,3,4,5)$
\item {\bf action:} $\Z/(2)$ acts by $(+,-,+,-,-,-)$
\item $\Bt=\{\Big({1\over 6}(1,1,5)\Big)_3,\Big({1\over 10}(1,3,7)\Big)_5\}=(\Bt 2.7)$ 
\item $\B \backslash \Bt=\{{1\over 2}(1,1,1),{1\over 3}(1,1,2)\}$
\end{itemize}
\leftline{\bf No. 25}
\begin{itemize}
\item {\bf cover:} $Y_{8,10}\subset\P(1,2,3,4,4,5)$
\item {\bf action:} $\Z/(2)$ acts by $(+,-,-,+,-,-)$
\item $\Bt=\{\Big({1\over 1}(1,1,1)\Big)_1,2\times\Big({1\over 4}(1,1,3)\Big)_2,\Big({1\over 6}(1,1,5)\Big)_3\}=(\Bt 2.16)$ 
\item $\B \backslash \Bt=\{{1\over 2}(1,1,1),{1\over 4}(1,1,3)\}$
\end{itemize}
\leftline{\bf No. 26}
\begin{itemize}
\item {\bf cover:} $Y_{8,12}\subset\P(1,2,3,4,5,6)$
\item {\bf action:} $\Z/(2)$ acts by $(+,+,-,-,-,-)$
\item $\Bt=\{3\times\Big({1\over 1}(1,1,1)\Big)_1,\Big({1\over 10}(1,1,9)\Big)_5\}=(\Bt 2.11)$ 
\item $\B \backslash \Bt=\{2\times{1\over 2}(1,1,1),{1\over 3}(1,1,2)\}$
\end{itemize}
\leftline{\bf No. 27}
\begin{itemize}
\item {\bf cover:} $Y_{10,12}\subset\P(1,2,3,5,5,7)$
\item {\bf action:} $\Z/(2)$ acts by $(+,-,-,+,-,-)$
\item $\Bt=\{\Big({1\over 1}(1,1,1)\Big)_1,\Big({1\over 14}(1,5,9)\Big)_7\}=(\Bt 2.3)$ 
\item $\B \backslash \Bt=\{{1\over 5}(1,2,3)\}$
\end{itemize}
\leftline{\bf No. 28}
\begin{itemize}
\item {\bf cover:} $Y_{10,12}\subset\P(1,3,4,4,5,6)$
\item {\bf action:} $\Z/(2)$ acts by $(+,-,+,-,-,-)$
\item $\Bt=\{2\times\Big({1\over 1}(1,1,1)\Big)_1,\Big({1\over 4}(1,1,3)\Big)_2,\Big({1\over 8}(1,1,7)\Big)_4\}=(\Bt 2.13)$ 
\item $\B \backslash \Bt=\{{1\over 3}(1,1,2),{1\over 4}(1,1,3)\}$
\end{itemize}

\vskip 10pt

\begin{rem}
Second degree (in the sense of Remark \ref{2nd-degree}) of all equations is 0 except in cases:
\begin{itemize}
\item {\bf No 4d}, where the equation of first degree 3 must have second degree equal to $2\in \Z/(4)$.
\item {\bf No 5b}, one of the equations must have second degree equal to $2\in \Z/(4)$.
\item {\bf No 6c}, one of the equations must have second degree equal to $2\in \Z/(4)$.
\item {\bf No 6d}, one of the equations must have second degree equal to $2\in \Z/(4)$.
\end{itemize}
\end{rem}

\begin{rem}\label{ej2}
Observe that case No. 13 is not a quotient from a codimension two Fletcher's example in {\bf [F]}. The Hilbert Series corresponds to an example in Reid's codimension one list: $Y_8\subset\P(1,1,2,2,3)$. However,  applying step (6) of the algorithm in Remark \ref{algor}, when we multiply in the Hilbert series of $X$ by $(1-t),(1-et),(1-t^2),(1-et^2)$ and $(1-et^3)$, we do not get $1-t^8$ or $1-et^8$ as expected. Instead we get 
$$1-t^4+et^4+\mathrm{higher\ degree}$$
This means that we need a generator of bidegree $(4,1)$ and an equation of bidegree $(4,0)$. Therefore, the new generator is not in the equation, and then the equation is not a linear cone, so our threefold cannot be projected ``quasismoothly" to a $\P(1,1,2,2,3)$ (for more details, see{\bf [F]}). It is a ``special" case described by Brown in {\bf [B]}. The rest of cases in Table 3 come all from Fletcher's list.

The same goes for No. 3 (whose Hilbert series corresponds to $Y_4\subset\P^4$).
\end{rem}

\vskip 10pt
Now we  list codimension 3 examples. Surprisingly enough, it contains only four cases, in contrast with the initial codimension-3 Fano threefolds, which are the intersection of 3 quadrics in $\P^6$ and the 69 examples in Alt{\i}nok's list. Observe also that from {\bf [A]} or {\bf [CR]} we know that last 69 are defined by the Pfaffians of a $5\times5$ skew-symmetric matrix. This is useful to search for equations. For the only case of Fano--Enriques which comes from a Fano in Alt{\i}nok's list ({\bf No.2} in Table 4 below), we  change the notation and write for the cover $Y_{d_1,...,d_5}$ to give the degrees of the equations of the pfaffians. Also, we add the specifications about the second degree of the equations of Remark \ref{2nd-degree}.

\vskip 10pt
\centerline{\bf Table 4. Fano--Enriques threefolds from codimension 3 Fano threefolds:}
\vskip 10pt
\leftline{\bf No. 1a}
\begin{itemize}
\item {\bf cover:} $Y_{2,2,2}\subset\P(1,1,1,1,1,1,1)$
\item {\bf action:} $\Z/(2)$ acts by $(+,+,+,-,-,-,-)$
\item $\Bt=\{8\times\Big({1\over 2}(1,1,1)\Big)_1\}=(\Bt 2.20)$ 
\item $\B \backslash \Bt=\emptyset$
\item Second degree of the equations must be 0
\end{itemize}
\leftline{\bf No. 1b}
\begin{itemize}
\item {\bf cover:} $Y_{2,2,2}\subset\P(1,1,1,1,1,1,1)$
\item {\bf action:} $\Z/(4)$ acts by $(+,+,i,i,-,-i,-i)$
\item $\Bt=\{2\times\Big({1\over 2}(1,1,1)\Big)_1,2\times\Big({1\over 4}(1,1,3)\Big)_1,2\times\Big({1\over 4}(1,1,3)\Big)_3\}=(\Bt 4.5)$ 
\item $\B \backslash \Bt=\emptyset$
\item Second degree of the equations must be 0, 0 and 2
\end{itemize}
\leftline{\bf No. 1c}
\begin{itemize}
\item {\bf cover:} $Y_{2,2,2}\subset\P(1,1,1,1,1,1,1)$
\item {\bf action:} $\Z/(8)$ acts by $(1,\epsilon,\epsilon^2,\epsilon^3,\epsilon^4,\epsilon^5,\epsilon^7)$
\item $\Bt=\{2\times\Big({1\over 2}(1,1,1)\Big)_1,\Big({1\over 4}(1,1,3)\Big)_1,\Big({1\over 8}(1,3,5)\Big)_3,\Big({1\over 8}(1,3,5)\Big)_7\}=(\Bt 8.1)$ 
\item $\B \backslash \Bt=\emptyset$
\item Second degree of the equations must be 0, 2 and 4
\end{itemize}
\leftline{\bf No. 2}
\begin{itemize}
\item {\bf cover:} $Y_{3,3,3,3,4}\subset\P(1,1,1,2,2,2,2)$
\item {\bf action:} $\Z/(5)$ acts by $(1,\epsilon,\epsilon^4,\epsilon,\epsilon^2,\epsilon^3,\epsilon^4)$
\item $\Bt=\{\Big({1\over 5}(1,1,4)\Big)_1,\Big({1\over 5}(1,1,4)\Big)_4,\Big({1\over 5}(1,2,3)\Big)_1,\Big({1\over 5}(1,2,3)\Big)_4\}=(\Bt 5.3)$ 
\item $\B \backslash \Bt=\{ {1\over 2}(1,1,1)\}$
\item Second degree of the equations must be 0, 1, 2, 3 and 4
\end{itemize}

\begin{rem}
Numbers 4 and 18 in Alt{\i}nok's list fit  the basket and $-K_Y^3$ restrictions when applying steps $(1)-(5)$ in the algorithm in Remark \ref{algor}. However, in analogy with Remark \ref{ej2}, their respective Hilbert series require a degree 4 relation and a degree 4 generator (with different second degrees), so one expects they are codimension 4 special cases.
\end{rem}
\vskip 10pt



\vskip 20pt
\centerline{\bf References}

\vskip 10pt

{\bf [A]} S. Alt{\i}nok. {\it Graded rings corresponding to polarised K3 surfaces and Fano threefolds}, Univ. of Warwick PhD thesis, Sep 1998, 93 + vii pp., get from www.maths.warwick.ac.uk/\~{}miles/doctors/Selma 

{\bf [ABR]} S. Alt{\i}nok, G. Brown and M. Reid. {\it Fano threefolds, K3 surfaces and graded rings.} In {\it Topology and geometry: commemorating SISTAG} (National Univ. of Singapore, 2001), Ed. A. J. Berrick and others, Contemp. Math. 314, AMS, 2002, pp. 25-53.

{\bf [B]} G. Brown. {\it Finding special K3 surfaces with Magma.} Preprint

{\bf [CR]} A. Corti and M. Reid. {\it Weighted grassmannians.} In {\it Algebraic geometry} (Genova, Sep 2001), In memory of Paolo Francia, M. Beltrametti and F. Catanese Eds., de Gruyter, Berlin, 2002, 141-163.

{\bf [D]} I. Dolgachev. {\it Weighted projective varieties.} In {\it Group actions and vector fields} (Vancouver B.C. 1981) LNM956, pp. 34-71.

{\bf [F]} A.R. Iano-Fletcher- {\it Working with weighted complete intersections.} In {\it Explicit birational geometry of threefolds}, CUP 2000, pp 101-173.

{\bf [GRDW]} G. Brown, Graded rings database webpage, see www.maths.warwick.ac.uk/grdb/

{\bf [M]} Magma (john Cannon's computer algebra system): W.Bosma, J. Cannon and C. Playoust, The Magma algebra system I: The user language, J. Symb. Comp 24 (1997) 235-265. See also www.maths.usyd.edu.au:8000/u/magma

\vskip 10pt

\end{document}